\documentclass[12pt]{amsart}
\usepackage{amssymb}
\usepackage{amsmath, amscd}
\usepackage{latexsym}
\usepackage[dvips]{graphicx}
\usepackage{lineno}
\usepackage[dvips]{color}
\numberwithin{equation}{section}
\newtheorem{theorem}{Theorem}[section]

\newtheorem{lemma}[theorem]{Lemma}

\newtheorem{definition}[theorem]{Definition}

\newtheorem{question}[theorem]{Question}

\begin{document}
\title{ \bf An upper bound on the Wiener index of a $k$-connected graph}

\author{Zhongyuan Che}
\address{Department of Mathematics, Penn State University\\
Beaver Campus,  Monaca, PA 15061, U.S.A.}
\email{zxc10@psu.edu}
\author{Karen L. Collins}
\address{Department of Mathematics and Computer Science\\
Wesleyan University, Middletown, CT 06459,  U.S.A.}
\email{kcollins@wesleyan.edu}

\begin{abstract}  
The Wiener index of a connected graph 
is the summation of all distances between unordered pairs of vertices of the graph.  
In this paper, we give an upper bound on the Wiener index
of a $k$-connected graph $G$ of order $n$ 
for integers $n-1>k \ge 1$:
 \[W(G) \le \frac{1}{4} n \lfloor \frac{n+k-2}{k} \rfloor (2n+k-2-k\lfloor \frac{n+k-2}{k} \rfloor).\]
Moreover, we show that 
this upper bound is sharp when $k \ge 2$ is even,  and can be obtained by the Wiener index of 
Harary graph $H_{k,n}$. 

\vskip 0.2in \noindent {\emph{keywords}}:
Harary graph,  $k$-connected graph, Wiener index 
\end{abstract}
\maketitle

\section{Introduction and Results}

The {\it Wiener index} $W(G)$ of a graph $G$ 
was first introduced by Wiener in 1947 for applications in chemistry \cite{W47} on studying the boiling points of paraffins. 
It is the summation of all distances between unordered pairs of vertices of $G$.
The concept of Wiener index has been studied under different names such as 
the {\it total status} by Harary \cite{H59}, the {\it total distance}  by Entringer et al. \cite{EJS76},  
and the {\it transmission} by Plesn\'ik \cite{P84} for various applications to topics including chemistry, 
communication, facility location, and cryptology. 
Due to its strong connection to chemistry, 
Wiener indices of trees \cite{DEG01} and Wiener indices of hexagonal systems \cite{{DGKZ02}} 
have been studied intensively.  
After more than 60 years from its birth, 
the research on Wiener index is still very lively. 
For instance, here is a list of some recent work: 
characterizations of trees with specified order and degree sequence that maximize the Wiener index \cite{SW12},  
the maximum Wiener index of unicyclic graphs with fixed maximum degree \cite{DZ12},
inverse Wiener index problems that search for trees with a given Wiener index \cite{FLS12},
Wiener indices of iterated line graphs of trees \cite{KPS12, KPS13, KPS13Ars}, 
Wiener indices of random trees \cite{W12}, and
Wiener index versus Szeged index in networks \cite{KN13}.

For most general classes of graphs,  there is no closed formula to calculate their Wiener indices, 
not even a recursive formula.
Finding bounds on Wiener indices of a general class of graphs has attracted many researchers' interest. 
Entringer et al. \cite{EJS76} showed that for any connected graph with a given order, 
the Wiener index is minimized by that of a complete graph and maximized by that of a path,
and the Wiener index of a tree with a given order attains the minimum value when it is a star 
and the maximum value when it is a path.
Walikar et al. \cite{WSR04} gave some bounds 
on the Wiener index of a graph in terms of the graph order 
together with one or two other graph parameters such as size, radius, diameter,  
connectivity, independent number, and chromatic number. 
Balakrishnan et al. \cite{BSL10} gave
a sharp lower bound on the Wiener index of an arbitrary graph in terms of three graph parameters altogether: 
order, size and diameter. 
For all integers $n-1>k\ge 1$, Gutman et al. \cite{GZ06} gave a sharp lower bound
on the Wiener index of a $k$-connected graph
(resp.,  a $k$-edge-connected graph) of order $n$.
They raised the question on finding an upper bound on the Wiener index of a $k$-connected graph
(resp., a $k$-edge-connected graph) of order $n$, and pointed out that it seems 
much more difficult. 

Motivated from their work, we provide an upper bound 
on the Wiener index of a $k$-connected graph of order $n$.  Our main 
theorem, proved in Section~4, is the following.

\vspace{.2in}
\noindent {\bf Theorem 4.3} Let $n-1 > k\geq1$ and $G$ be a 
$k$-connected graph of order $n$.  Then
  \[W(G) \le \frac{1}{4} n \lfloor \frac{n+k-2}{k} \rfloor 
(2n+k-2-k\lfloor \frac{n+k-2}{k} \rfloor).\]

We show further that this upper bound is sharp when $k \ge 2$ 
is even. 
It is easily seen that adding an edge in a connected graph
decreases the Wiener index.
Thus, when we look for an upper bound on the Wiener index of a 
specific class of graphs of a given order, it is natural to consider those 
graphs with the minimum number of edges.
For example, a path of order $n$ has the minimum 
number of edges in the class of $1$-connected graphs of order $n$, and also has 
the maximum Wiener index in this class of graphs.  
We recall that in Section ~2 the well-known Harary graph, $H_{k,n}$,
where integers $n>k \ge 2$, has the minimum 
number of edges in the class of $k$-connected graphs of order $n$. 
In Section~3, we prove that if $k \ge 2$ is even, then the Wiener index of $H_{k,n}$ is equal to the 
maximum value in Theorem 4.3, and so the upper bound given 
by the theorem is tight for all positive even integers $k$.

The outline of the paper is as follows.  In Section 2, we give 
definitions and preliminaries on Wiener index and Harary graph $H_{k,n}$.
We recall that the {\it status of a vertex} \cite{H59} is the summation of 
distances between the vertex and all other vertices of the graph. 
In Section 3, we calculate the status of each vertex in $H_{k,n}$.
We obtain the Wiener index of $H_{k,n}$ as half of the summation of all vertex statuses of $H_{k,n}$.
In Section 4,  we provide an upper bound on the status of 
any vertex in a $k$-connected graph of order $n$ for integers $n-1>k\ge 1$
and apply it to prove Theorem~4.3. 

\section{Preliminaries}
All graphs in this paper are finite, simple and connected.
The cardinality of a set $S$  is denoted by $|S|$.
The vertex set of a graph $G$ is denoted by $V(G)$, 
and its cardinality is called the {\it order} of $G$.
A graph $G$ is called {\it connected} if any two vertices are joined by a path in $G$. 
The {\it distance} $d_{G}(x,y)$ between two vertices $x,y$ of $G$ is the length of a shortest path 
joining $x$ and $y$ in $G$.
The set of all vertices with distance $i$ from a vertex $x$ in $G$ is
denoted by $N(x,i)$. 
In particular, when $i=1$,  $N(x,1)$ is the set of neighborhoods of $x$ and its cardinality is called the {\it degree} of $x$. 
The {\it eccentricity of a vertex} $x$ of $G$, denoted by ecc$(x)$, 
is the maximum distance between $x$ and another vertex of $G$.
The {\it diameter of a graph} $G$ is the maximum vertex eccentricity of $G$ and is denoted by diam$(G)$.
The {\it status of a vertex} $x$ of $G$, denoted by $W(x,G)$, is the summation of all distances between
$x$ and other vertices of $G$, that is, \[W(x, G)= \sum_{y \in V(G)} d_G(x, y).\]
The {\it Wiener index} of a graph $G$, denoted by $W(G)$, is the summation of all distances between 
unordered pairs of vertices of $G$.  
So $W(G)$ can be written as:
\[W(G)= \sum\limits_{\{x,y\} \subseteq V(G)} d_{G}(x, y)=\frac{1}{2}\sum\limits_{x \in V(G)} W(x, G).\]

The following property of the vertex status comes immediately from its definition, 
and we use it as a key formula to find an upper bound on the Wiener index of a $k$-connected graph of order $n$.

\begin{lemma}\rm\cite{BF90}\label{L:W(x,G)N(x,i)}
Let $x$ be a vertex of a connected graph $G$ with the eccentricity ecc$(x)$. 
Assume that $N(x,i)$ is the set of all vertices with distance 
$i$ from $x$ in $G$, where $1 \le i \le \mbox{ecc}(x)$.  Then the status of the vertex $x$ is  
\[W(x, G)=\sum\limits_{i=1}^{ecc(x)} i \cdot |N(x,i)|.\]
\end{lemma}
 
A connected graph $G$ is called {\it $k$-connected} if the removal of 
any $k-1$ vertices of $G$ does not result a disconnected graph 
or a $1$-vertex graph. By convention, any complete graph $K_{n}$ 
is $(n-1)$-connected.
A {\it vertex cut}  of a connected graph $G$ is a set of vertices 
whose removal disconnects $G$.
Therefore, if $G$ is a $k$-connected graph with some nonadjacent vertices,
then any vertex cut of $G$ has at least $k$ vertices.
It is well known \cite{W00} that for integers $n>k \ge 2$,
any $k$-connected graph of order $n$ has 
at least $\lceil \frac{kn}{2} \rceil$ edges, 
and Harary graph $H_{k,n}$ is a $k$-connected graph of order $n$ 
that has exactly  $\lceil \frac{kn}{2} \rceil$ edges.    

\begin{definition} \rm \cite{H62}
Let $n> k \ge 2$. Place $n$ vertices $0, 1, 2, \cdots, n-1$ around a circle in the clockwise direction 
and equally spaced. The construction of Harary graph $H_{k,n}$ 
depends on the parity of $k$ and $n$ and falls into three cases:

Case 1. $k$ is even. Construct $H_{k,n}$ by making each vertex adjacent to the nearest $\frac{k}{2}$ vertices in
each direction around the circle.

Case 2. $k$ is odd and $n$ is even. Construct $H_{k,n}$ by making each vertex adjacent to the nearest $\frac{k-1}{2}$
vertices in each direction and to the diametrically opposite vertex.

Case 3. both $k$ and $n$ are odd. Construct $H_{k,n}$ from $H_{k-1,n}$ by adding the edge between vertices $i$ and 
$i+\frac{n-1}{2}$ for $0 \le i \le \frac{n-1}{2}$.
\end{definition}

By the definition of Harary graph, we can see that if at least one of $k$ and $n$ is even, 
then $H_{k,n}$ is vertex-transitive, and each vertex has degree $k$.
On the other hand, if both $k$ and $n$ are odd, then $H_{k,n}$ is not vertex-transitive, 
it has exactly one vertex, namely vertex ``$\frac{n-1}{2}$'',  with degree $k+1$ and all other vertices have degree $k$.

\section{Harary graphs $H_{k,n}$} 

In this section, we calculate diam$(H_{k,n})$,  $W(x,H_{k,n})$ and $W(H_{k,n})$ for Harary graph
$H_{k,n}$ where $n-1>k \ge 2$. We exclude the case when $n=k+1$, that is, when
$H_{k,n}$ is a complete graph. 
Though the Wiener index of Harary graph $H_{k,n}$ attains
the upper bound given by Theorem \ref{T:W(G)} only when $k \ge 2$ is even, 
we provide $W(H_{k,n})$ for all $n-1>k \ge 2$ for the completeness of studying this topic.

Let $n-1>k \ge 2$ where $k$ is even. Then $H_{k,n}$ is vertex-transitive.
So $\mbox{ecc}(x)=\mbox{diam}(H_{k,n}) >1$ for any vertex $x$ of $H_{k,n}$.
We claim that $|N(x,i)| = k$ for $1 \le i \le \mbox{ecc}(x)-1$ as follows.
It is easy to see that $|N(x,1)|=k$ since $x$ is adjacent to the nearest $\frac{k}{2}$ vertices in
both directions around the circle.
Each vertex in $N(x,1)$ is adjacent to their nearest $\frac{k}{2}$ vertices in
both directions  around the circle. Then the neighborhoods of vertices in $N(x,1)$
contribute $k$ vertices to $N(x,2)$, which are the next nearest $\frac{k}{2}$ vertices to $x$ in
both directions around the circle.
Continue this way, we can see that $|N(x,i)| = k$ for $1 \le i \le \mbox{ecc}(x)-1$.
Finally, when $i=\mbox{ecc}(x)=t$, the neighborhoods of vertices in $N(x,t-1)$
contribute at least $1$ vertex and at most $k$ vertices to $N(x,t)$, 
and so $1 \le |N(x,t)| \le k$.

\begin{lemma}\label{L:H(k,n)kEven}\rm
Let $n-1>k \ge 2$ where $k$ is even. Then
\begin{eqnarray*}
\mbox{diam}(H_{k,n}) &=& D=\lfloor \frac{n+k-2}{k} \rfloor,\\
W(x,H_{k,n})&=& \frac{1}{2}D (2n+k-2-kD), \\
W(H_{k,n})&=& \frac{1}{4} n D (2n+k-2-kD).
\end{eqnarray*}
\end{lemma}
\proof  If $k \ge 2$ is even, then $H_{k,n}$ is vertex-transitive.
Assume that each vertex $x$ of $H_{k,n}$ has eccentricity $\mbox{ecc}(x)=t \ge 2$.
Note that $n=|V(H_{k,n})| =\sum\limits_{i=0}^{t} |N(x,i)|$, where
$|N(x,0)|=1$ and $|N(x,i)|=k$ for $1 \le i \le t-1$,
and $1 \le |N(x,t)| \le k$.  
Then $1+k(t-1)+1 \le n$ and $n \le 1+k(t-1)+ k$. 
It follows that 
$\lceil \frac{n-1}{k}\rceil \le t \le \lfloor \frac{n+k-2}{k} \rfloor.$

Let $n=mk+i$ where $m \ge 1$ and $0 \le i \le k-1$. 
Then 
\begin{eqnarray*}
\lceil \frac{n-1}{k}\rceil &=&\lceil m+ \frac{i-1}{k}\rceil 
=\left\{\begin{array}{cc}
m, & \mbox{if $i=0,1$},\\
m+1, & \mbox{if $2 \le i \le k-1$}.
\end{array}
\right.\\
\lfloor \frac{n+k-2}{k} \rfloor &= &\lfloor m+1 + \frac{i-2}{k} \rfloor 
=\left\{\begin{array}{cc}
m, & \mbox{if $i=0, 1$},\\
m+1, & \mbox{if $2 \le i \le k-1$}.
\end{array}
\right.
\end{eqnarray*}

Therefore, $t=\lceil \frac{n-1}{k}\rceil = \lfloor \frac{n+k-2}{k} \rfloor$, and
$\mbox{diam}(H_{k,n})=t=\lfloor \frac{n+k-2}{k} \rfloor.$

We now calculate $W(x,H_{k,n})$.  
By Lemma \ref{L:W(x,G)N(x,i)}, 
\begin{eqnarray*}
W(x,H_{k,n}) &=& \sum\limits_{i=1}^{t} i \cdot |N(x,i)|= \sum\limits_{i=1}^{t} (t-(t-i)) \cdot |N(x,i)| \\
&=& t \sum\limits_{i=1}^{t} |N(x,i)| - \sum\limits_{i=1}^{t-1} (t-i) |N(x,i)|.
\end{eqnarray*}

Since $\sum\limits_{i=1}^{t} |N(x,i)|=n-1$ and $|N(x,i)| = k$ for $1 \le i \le t-1$,
\begin{eqnarray*}
W(x,H_{k,n})&=& (n-1) t - k \sum\limits_{i=1}^{t-1}(t-i)= (n-1) t - k \sum\limits_{i=1}^{t-1}i \\
&=& (n-1)t -k {t \choose 2} \\
&=& -\frac{k}{2} t^{2} +\frac{2n+k-2}{2} t\\
&=& \frac{1}{2}\lfloor \frac{n+k-2}{k} \rfloor (2n+k-2-k\lfloor \frac{n+k-2}{k} \rfloor). 
\end{eqnarray*}

Recall that $H_{k,n}$ is vertex-transitive. Then each vertex has the same vertex status $W(x, H_{k,n})$. 
Hence, 
\begin{eqnarray*}
W(H_{k,n})&=& \frac{1}{2}\sum_{x \in V(H_{k,n})} W(x, H_{k,n})=\frac{1}{2}n W(x, H_{k,n})\\
&=&\frac{1}{4} n \lfloor \frac{n+k-2}{k} \rfloor (2n+k-2-k \lfloor \frac{n+k-2}{k} \rfloor).
\end{eqnarray*}
\qed\\

Let $n-1>k \ge 3$ where $k$ is odd and $n$ is even. Then $H_{k,n}$ is vertex-transitive. 
So $\mbox{ecc}(x)=\mbox{diam}(H_{k,n})>1$ for any vertex $x$ of $H_{k,n}$.
It is clear that $|N(x,1)|=k$ for any vertex $x$ of $H_{k,n}$.
If $\mbox{diam}(H_{k,n})=2$, then $|N(x,2)|=n-k-1$.
Suppose that $\mbox{diam}(H_{k,n})=t>2$. 
Then $|N(x,i)|=2(k-1)>k$ for $2 \le i \le t-1$ as follows.
By definition, each vertex $x$ is adjacent to the nearest $\frac{k-1}{2}$
vertices in both directions around the circle, and to its diametrically opposite vertex $x'=x+\frac{n}{2}$. 
So $|N(x,1)|=k$. 
Then the neighborhoods of $x'$ contribute $k-1$ vertices to $N(x,2)$ 
which are the nearest $\frac{k-1}{2}$ vertices to $x'$ in both directions around the circle. 
Moreover, the neighborhoods of vertices in $N(x,1)\setminus \{x'\}$ 
will contribute $k-1$ more vertices to $N(x,2)$,  
which are the next nearest  $\frac{k-1}{2}$ vertices to $x$ in both directions around the circle.
Therefore, $|N(x,2)|=2(k-1)$. 
Continue this way, we can see that $|N(x,i)| =2(k-1)$ for $2 \le i \le t-1$.
Finally, when $i=\mbox{ecc}(x)=t$,  
$|N(x,t)|$ is even since $|N(x,t)|=n-[1+k+\sum\limits_{i=2}^{t-1} 2(k-1)]$ 
where $k$ is odd and $n$ is even.
Hence,  $2 \le |N(x,t)| \le 2(k-1)$. 

\begin{lemma}\rm \label{L:H(k,n)kOddnEven}
Let $n-1>k \ge 3$ where $k$ is odd and $n$ is even. 

(i) If $k+1 < n \le 3k-1$, then $\mbox{diam}(H_{k,n}) = 2$, and
\begin{eqnarray*}
W(x,H_{k,n})&=& 2n-k-2, \\
W(H_{k,n})&=& \frac{1}{2} n (2n-k-2).
\end{eqnarray*}

(ii) If $n \ge 3k+1$, then  $\mbox{diam}(H_{k,n}) \ge 3$, and
\begin{eqnarray*}
\mbox{diam}(H_{k,n})&=&D=\lfloor \frac{n-k-3}{2(k-1)} \rfloor+2,\\
W(x, H_{k,n}) &=& \frac{1}{2} D (2n+4k-8 -2(k-1)D)-(k-2),\\
W(H_{k,n}) &=& \frac{1}{4} n D (2n+4k-8 -2(k-1)D)- \frac{1}{2}n(k-2).\\
\end{eqnarray*}
\end{lemma}

\proof Let $n-1>k \ge 3$ where $k$ is odd and $n$ is even. Then
$H_{k,n}$ is a vertex-transitive graph. 

If $k+1 < n \le 3k-1$, then by the definition of $H_{k,n}$, 
it is easy to check that each vertex $x$ of $H_{k,n}$ has eccentricity $\mbox{ecc}(x)=t=2$,
and so $\mbox{diam}(H_{k,n}) = 2$.
Moreover, $|N(x,1)|=k$ and $|N(x,2)|=n-k-1$. It follows that
$W(x,H_{k,n})= 1 \cdot k +2(n-k-1)=2n-k-2$, and 
$W(H_{k,n}) = \frac{1}{2} n (2n-k-2)$.

Note that $n \neq 3k$ since $k$ is odd and $n$ is even.
From now on, we assume that $n \ge 3k+1$. By definition of $H_{k,n}$, it is easy to check that 
each vertex $x$ of $H_{k,n}$ has $\mbox{ecc}(x)=t \ge 3$.
We have known that $|N(x,0)|=1$, $|N(x,1)|=k$, $|N(x,i)|=2(k-1)$ for $2 \le i \le t-1$,
and $2 \le |N(x,t)| \le 2(k-1)$.  
By $n=|V(H_{k,n})| =\sum\limits_{i=0}^{t} |N(x,i)|$,
we have that $1+k+2(k-1)(t-2)+2 \le n$ and $n \le 1+k+2(k-1)(t-2)+ 2(k-1)$. 
It follows that 
$\lceil \frac{n-k-1}{2(k-1)} \rceil +1 \le t \le \lfloor \frac{n-k-3}{2(k-1)} \rfloor+2.$
Note that $n-k-1$ is even since $k$ is odd and $n$ is even. 
We can write $\frac{n-k-1}{2} = m(k-1)+i$ for some integers $m \ge 0$ and $0 \le i \le k-2$.
Then 
\begin{eqnarray*}
\lceil \frac{n-k-1}{2(k-1)} \rceil +1&=&\lceil m+\frac{i}{k-1} \rceil +1
=\left\{\begin{array}{cc}
m+1, & \mbox{if $i=0$},\\
m+2, & \mbox{if $1 \le i \le k-2$}.
\end{array}
\right.\\
 \lfloor \frac{n-k-3}{2(k-1)} \rfloor+2 &= &\lfloor m+\frac{i-1}{k-1} \rfloor +2
 =\left\{\begin{array}{cc}
m+1, & \mbox{if $i=0$},\\
m+2, & \mbox{if $1 \le i \le k-2$}.
\end{array}
\right.
\end{eqnarray*}

Therefore, $t=\lceil \frac{n-k-1}{2(k-1)} \rceil +1= \lfloor \frac{n-k-3}{2(k-1)} \rfloor+2$, and
\begin{eqnarray*}
\mbox{diam}(H_{k,n})&=&t= \lfloor \frac{n-k-3}{2(k-1)} \rfloor+2.
\end{eqnarray*}

We then calculate $W(x,H_{k,n})$. By Lemma \ref{L:W(x,G)N(x,i)},
\begin{eqnarray*}
W(x,H_{k,n}) &=& \sum\limits_{i=1}^{t} i \cdot |N(x,i)|, \mbox{ where $t=\mbox{ecc}(x)$}\\
&=& t \sum\limits_{i=1}^{t} |N(x,i)| - \sum\limits_{i=1}^{t-1} (t-i) |N(x,i)|\\
&=& (n-1) t - k (t-1)- 2(k-1)\sum\limits_{i=2}^{t-1}(t-i)\\
&=& (n-1) t - k(t-1)-2(k-1) \sum\limits_{i=1}^{t-2}i \\
&=& (n-1)t -k (t-1) - 2(k-1) {t-1 \choose 2} \\
&=& -(k-1) t^{2} + (n+2k-4)t + (-k+2).
\end{eqnarray*}

Recall that $t=\mbox{diam}(H_{k,n})=\lfloor \frac{n-k-3}{2(k-1)} \rfloor+2$.
For clarity, we denote $D=t=\mbox{diam}(H_{k,n})$
and rewrite $W(x,H_{k,n})$ using $D$.
\begin{eqnarray*}
W(x,H_{k,n}) &=&-(k-1) D^{2} + (n+2k-4)D + (-k+2)\\
&=&\frac{1}{2} D (2n+4k-8 -2(k-1)D)-(k-2).
\end{eqnarray*}

Since $H_{k,n}$ is vertex-transitive, each vertex has the same vertex status $W(x, H_{k,n})$. 
Hence, 
\begin{eqnarray*}
W(H_{k,n}) &=& \frac{1}{2} \sum_{x \in V(H_{k,n})} W(x, H_{k,n})=\frac{1}{2} n W(x, H_{k,n})\\
&=&\frac{1}{4} n D (2n+4k-8 -2(k-1)D)- \frac{1}{2}n(k-2).
\end{eqnarray*}
\qed\\

Let $n-1>k \ge 3$ where both $k$ and $n$ are odd integers. 
Then $H_{k,n}$ is not vertex-transitive. 
It has exactly one vertex of degree $k+1$ and all other vertices have degree $k$.
By definition. $H_{k,n}$ is constructed from $H_{k-1,n}$ by adding an edge between vertices $x$ and 
$x+\frac{n-1}{2}$ for $0 \le x \le \frac{n-1}{2}$. 

If $x=\frac{n-1}{2}$, then the vertex $x$ is adjacent to the nearest $\frac{k-1}{2}$
vertices in both directions around the circle, and to two diametrically opposite vertices $x'=0$ and $x''=n-1$.
So $|N(x,1)|=k+1$.  If $\mbox{ecc}(x) \le 2$, then $|N(x,2)|=n-k-2$, which is $0$ when $n=k+2$.
Suppose that $\mbox{ecc}(x)=t>2$.  Then
similarly to the argument for the case when $k$ is odd and $n$ is even, 
we can see that $|N(x,i)| =2(k-1)$ for $2 \le i \le t-1$.
Finally, when $i=\mbox{ecc}(x)=t$, $|N(x,t)|$ is even since $|N(x,t)|=n-[k+2+\sum\limits_{i=2}^{t-1} 2(k-1)]$
where both $k$ and $n$ are odd integers. 
Hence,  $2 \le |N(x,t)| \le 2(k-1)$. 

If $x \neq \frac{n-1}{2}$, then the vertex $x$ is adjacent to the nearest $\frac{k-1}{2}$
vertices in both directions around the circle, and to one diametrically opposite vertex $x'=x+\frac{n-1}{2}$.
So $|N(x,1)|=k$ and $\mbox{ecc}(x) \ge 2$ since $n>k+1$.  
If $\mbox{ecc}(x)=2$, then $|N(x,2)|=n-k-1$.
Suppose that $\mbox{ecc}(x)=t>2$.  Then
similarly to the argument for the case when $k$ is odd and $n$ is even, 
we can see that $|N(x,i)| =2(k-1)$ for $2 \le i \le t-1$.
Finally, when $i=\mbox{ecc}(x)=t$, $|N(x,t)|$ is odd since $|N(x,t)|=n-[k+1+\sum\limits_{i=2}^{t-1} 2(k-1)]$
where both $k$ and $n$ are odd integers. 
Hence, $1 \le |N(x,t)| \le 2k-3$.

\begin{lemma}\rm \label{L:H(k,n)kOddnOdd}
Let $n-1>k \ge 3$ where both $k$ and $n$ are odd integers. 
Assume that $z$ is the vertex of degree $k+1$
and $x$ is a vertex of degree $k$ in $H_{k,n}$.

(i) If $k+1 < n \le 3k-2$, then $\mbox{diam}(H_{k,n})=2$, and 
\begin{eqnarray*}
W(z, H_{k,n})&=&2n-k-3,\\
W(x,H_{k,n})&=& 2n-k-2, \\
W(H_{k,n})&=& \frac{1}{2} n (2n-k-2)-\frac{1}{2}.
\end{eqnarray*}

(ii) If $n=3k$, then $\mbox{diam}(H_{k,n})=3$, and
\begin{eqnarray*}
W(z, H_{k,n})&=&2n-k-3,\\
W(x,H_{k,n})&=& 2n-k-1, \\
W(H_{k,n})&=& \frac{1}{2} n (2n-k-1)-1.
\end{eqnarray*}

(iii) If $n \ge 3k+2$, then $\mbox{diam}(H_{k,n}) \ge 3$, and
\begin{eqnarray*}
\mbox{diam}(H_{k,n})&=&D=\lfloor \frac{n-k-2}{2(k-1)} \rfloor +2,\\
W(z, H_{k,n})&=&\frac{1}{2}D(2n+4k-10-2(k-1)D) -(k-3),\\
W(x,H_{k,n})&=&\frac{1}{2} D (2n+4k-8-2(k-1)D) -(k-2),\\
W(H_{k,n})&=&\frac{1}{4}n D(2n+4k-8 - 2(k-1)D)-\frac{1}{2}(n(k-2)+D-1).
\end{eqnarray*}
\end{lemma}

\proof Let $n-1>k \ge 3$ where both $k$ and $n$ are odd integers. 
Note that $H_{k,n}$ is not vertex-transitive 
since it has exactly one vertex of degree $k+1$ and all other vertices have degree $k$.
Let $z$ be the vertex of degree $k+1$ and $x$ be a vertex of degree $k$ in $H_{k,n}$.

If $k+1 < n  \le 3k-2$,  then by the definition of $H_{k,n}$, 
it is easy to check that the vertex of degree $k+1$ has eccentricity at most $2$,
and each  vertex of degree $k$ has eccentricity $2$.
Hence, $\mbox{diam}(H_{k,n})=2$.

Since $|N(z,1)|=k+1$ and $|N(z,2)|=n-k-2$, 
\[ W(z,H_{k,n}) = 1 \cdot (k+1) + 2(n-k-2)=2n-k-3.\]

Since $|N(x,1)|=k$ and $|N(x,2)|=n-k-1$,
\[ W(x,H_{k,n}) = 1 \cdot k + 2(n-k-1)=2n-k-2.\]

Note that $H_{k,n}$ has exactly one vertex with $W(z, H_{k,n})$ and 
$n-1$ vertices with $W(x,H_{k,n})$. 
Then 
\begin{eqnarray*}
W(H_{k,n})&=&\frac{1}{2} W(z,H_{k,n}) + \frac{1}{2}(n-1)W(x, H_{k,n})\\
&=&\frac{1}{2}( 2n-k-3) + \frac{1}{2} (n-1) (2n-k-2)\\
&=&\frac{1}{2}n (2n-k-2)-\frac{1}{2}.
\end{eqnarray*}

Note that $n \neq 3k-1, 3k+1$ since both $k$ and $n$ are odd integers.
If $n=3k$,  then by the definition of $H_{k,n}$, 
it is easy to check that  
the vertex of degree $k+1$ has eccentricity $2$,
and each vertex of degree $k$ has eccentricity $3$.
Hence, $\mbox{diam}(H_{k,n})=3$.

Since $|N(z,1)|=k+1$ and $|N(z,2)|=n-k-2$,
\[ W(z,H_{k,n}) = 1 \cdot (k+1) + 2\cdot(n-k-2)=2n-k-3.\]

Since $|N(x,1)|=k$, $|N(x,2)|=2(k-1)$ and $|N(x,3)|=n-3k+1$,
\[ W(x,H_{k,n}) = 1 \cdot k + 2\cdot 2(k-1) +3 \cdot (n-3k+1)=2n-k-1.\]

Note that $H_{k,n}$ has exactly one vertex with $W(z, H_{k,n})$ and 
$n-1$ vertices with $W(x,H_{k,n})$. 
Then 
\begin{eqnarray*}
W(H_{k,n})&=&\frac{1}{2} W(z,H_{k,n}) + \frac{1}{2}(n-1)W(x, H_{k,n})\\
&=&\frac{1}{2}( 2n-k-3) + \frac{1}{2} (n-1) (2n-k-1)\\
&=&\frac{1}{2}n (2n-k-1)-1.
\end{eqnarray*}

For the rest of the proof, we assume that $n \ge 3k+2$.  Then $\mbox{ecc}(z), \mbox{ecc}(x) \ge 3$.
We will first calculate $\mbox{diam}(H_{k,n})$ by computing $\mbox{ecc}(z)$ and $\mbox{ecc}(x)$ respectively.

Let $\mbox{ecc}(z)=\bar{t} \ge 3$. 
Then $|N(z,0)|=1$, $|N(z,1)|=k+1$, $|N(z,i)|=2(k-1)$ for $2 \le i \le \bar{t}-1$, 
and $2 \le |N(z,\bar{t})| \le 2(k-1)$.  
By the fact that $n=|V(H_{k,n})| =\sum\limits_{i=0}^{\bar{t}} |N(z,i)|$,
we can see that $1+k+1+2(k-1)(\bar{t}-2)+2 \le n$ and  $n \le 1+k+1+2(k-1)(\bar{t}-2)+ 2(k-1)$. 
It follows that 
\[\lceil \frac{n-k-2}{2(k-1)} \rceil +1 \le \bar{t} \le \lfloor \frac{n-k-4}{2(k-1)} \rfloor+2.\] 

Note that $n-k-2$ is even. Then $\frac{n-k-2}{2}=m(k-1)+i$ for some integers 
$m \ge 1$ and $0 \le i \le k-2$.  Then
\begin{eqnarray*}
\lceil \frac{n-k-2}{2(k-1)} \rceil +1&=&\lceil m+\frac{i}{k-1} \rceil +1\\
&=&\left\{
\begin{array}{cc}
m+1, & \mbox{if $i=0$}, \\
m+2, & \mbox{if $1 \le i \le k-2$}.
\end{array}
\right.\\
\lfloor \frac{n-k-4}{2(k-1)} \rfloor+2 &= &\lfloor m+\frac{i-1}{k-1} \rfloor +2\\
 &=&\left\{\begin{array}{cc}
m+1, & \mbox{if $i=0$}, \\
m+2, & \mbox{if $1 \le i \le k-2$}.
\end{array}
\right.
\end{eqnarray*}

Hence, $\mbox{ecc}(z)=\bar{t}=\left\{\begin{array}{cc}
m+1, & \mbox{if $i=0$}, \\
m+2, & \mbox{if $1 \le i \le k-2$}.
\end{array}
\right.$

Let $\mbox{ecc}(x)=t \ge 3$.
Then $|N(x,0)|=1$ and $|N(x,1)|=k$, $|N(x,i)|=2(k-1)$ for $2 \le i \le t-1$,
and $1 \le |N(x,t)| \le 2k-3$.  
By $n=|V(H_{k,n})| =\sum\limits_{i=0}^{t} |N(x,i)|$,
we have that $1+k+2(k-1)(t-2)+1 \le n$ and $n \le 1+k+2(k-1)(t-2)+ 2k-3$. It follows that 
\[\lceil \frac{n-k}{2(k-1)} \rceil +1 \le t \le \lfloor \frac{n-k-2}{2(k-1)} \rfloor+2.\] 

Recall that $\frac{n-k-2}{2}=m(k-1)+i$ where $m \ge 1$ and $0 \le i \le k-2$. Then
\begin{eqnarray*}
\lceil \frac{n-k}{2(k-1)} \rceil +1&=&\lceil m+\frac{i+1}{k-1} \rceil +1=m+2,\\
\lfloor \frac{n-k-2}{2(k-1)} \rfloor+2 &= &\lfloor m+\frac{i}{k-1} \rfloor +2=m+2.
 \end{eqnarray*}

Therefore,  $\mbox{ecc}(x)=t=m+2$. 

Now we have $m+1 \le \mbox{ecc}(z) \le m+2$ and $\mbox{ecc}(x)=m+2$  
where $m=\lfloor \frac{n-k-2}{2(k-1)} \rfloor \ge 1$. 
Therefore,
\[\mbox{diam}(H_{k,n})=D=m+2=\lfloor\frac{n-k-2}{2(k-1)}\rfloor+2.\]

We then calculate $W(z, H_{k,n})$ and $W(x,H_{k,n})$ respectively.

Similarly to the calculation in Lemma \ref{L:H(k,n)kOddnEven}, we can have
\begin{eqnarray*}
W(z,H_{k,n}) &=& \sum\limits_{i=1}^{\bar{t}} i \cdot |N(z,i)|, \mbox{where $\bar{t}=\mbox{ecc}(z) \ge 3$} \\
&=& -(k-1)\bar{t}^{2} + (n+2k-5) \bar{t} +(-k+3).
\end{eqnarray*}

Recall that $\bar{t}=m+1$ or $m+2$ based on $\frac{n-k-2}{2} \equiv 0 ($mod $k-1)$ or not. 
Then we distinguish two cases.

Case 1. If $\frac{n-k-2}{2} \not\equiv 0 ($mod $k-1)$, then $\bar{t}=m+2=D$.  
\begin{eqnarray*}
W(z,H_{k,n})&=&-(k-1)D^{2} + (n+2k-5) D +(-k+3)\\
&=&\frac{1}{2}D(2n+4k-10-2(k-1)D) -(k-3)
\end{eqnarray*}

Case 2. If $\frac{n-k-2}{2} \equiv 0 ($mod $k-1)$, then $\bar{t}=m+1=D-1$. 
\begin{eqnarray*}
W(z, H_{k,n})&=& -(k-1)(D-1)^{2} + (n+2k-5) (D-1) +(-k+3)\\
&=& -(k-1)D^{2} +(n+4k-7)D +(-n-4k+9)\\
&=& \frac{1}{2}D(2n+4k-10 -2(k-1)D) -(n+4k-9) +2D(k-1)\\
&=& \frac{1}{2}D(2n+4k-10 -2(k-1)D) -(k-3).
\end{eqnarray*}

The last equation can be obtained as follows.
Recall that $D=\lfloor\frac{n-k-2}{2(k-1)}\rfloor+2$. 
Hence, if $\frac{n-k-2}{2} \equiv 0 ($mod $k-1)$, then $D= \frac{n-k-2}{2(k-1)}+2$.
The term $2D(k-1)$ in the next to the last equation is 
$2D (k-1) = 2(\frac{n-k-2}{2(k-1)}+2)(k-1)=n+3k-6$. 
Simplify two terms $-(n+4k-9) +2D(k-1)$ in the next to the last equation
as $-(n+4k-9)+n+3k-6=-k+3$.
Then we get the last equation. 

Therefore, 
\[W(z, H_{k,n}) = \frac{1}{2}D(2n+4k-10 -2(k-1)D) -(k-3).\]

Similarly to the calculation in Lemma \ref{L:H(k,n)kOddnEven}, we can have
\begin{eqnarray*}
W(x,H_{k,n}) &=& \sum\limits_{i=1}^{t} i \cdot |N(x,i)|, \mbox{where $t=\mbox{ecc}(x)=D$}\\
&=&-(k-1)t^{2} + (n+2k-4) t + (-k+2) \\
&=& -(k-1) D^{2} + (n+2k-4)D+ (-k+2)\\
&=&\frac{1}{2} D (2n+4k-8-2(k-1)D) -(k-2).
\end{eqnarray*}

Finally, we will calculate $W(H_{k,n})$.
Note that $H_{k,n}$ has exactly one vertex with $W(z, H_{k,n})$ and 
$n-1$ vertices with $W(x,H_{k,n})$. 
Then 
\begin{eqnarray*}
W(H_{k,n})&=&\frac{1}{2} W(z,H_{k,n}) + \frac{1}{2}(n-1)W(x, H_{k,n})\\
&=&\frac{1}{2}( \frac{1}{2}D(2n+4k-10 -2(k-1)D) -(k-3))\\
&& + \frac{1}{2} (n-1) (\frac{1}{2} D (2n+4k-8-2(k-1)D) -(k-2))\\
&=&\frac{1}{4}n D(2n+4k-8 - 2(k-1)D)-\frac{1}{2}(n(k-2)+D-1).
\end{eqnarray*}
\qed\\

\section{$k$-connected graphs}

In this section, we give upper bounds on $W(x,G)$ and $W(G)$ 
for a $k$-connected graph $G$ of order $n$ where integers $n-1>k \ge 1$.
We exclude the case when $n=k+1$, that is, when $G$ is the complete graph $K_{n}$, whose
Wiener index is well-known as $\frac{1}{2} n^{2} - \frac{1}{2}n$.

The following result is given as Exercise 4.2.22 in \cite{W00}. 
We provide a proof for the completeness.
\begin{theorem}\rm \cite{W00} \label{T:D(G)}
Let  $n-1>k \ge 1$ and $G$ be a $k$-connected graph of order $n$. 
Then $\mbox{diam}(G) \le \lfloor \frac{n+k-2}{k} \rfloor$ and this bound is sharp when 
$k=1$ or $k \ge 2$ is even. 
\end{theorem}
\proof  Let $x$ be a vertex of $G$ with $\mbox{ecc}(x)=t$ where $t>1$.
Such a vertex exists since $G$ is not a complete graph.
Let $N(x,i)$ be the set of vertices with distance $i$ from $x$ in $G$ for $0 \le i \le t$.
Then $N(x,i)$ are pairwise disjoint for $0 \le i \le t$, and form a partition of the vertex set of $G$, that is,
$n=|V(G)| =\sum\limits_{i=0}^{t} |N(x,i)|$. 
It is clear that $|N(x,0)|=1$ and $|N(x,t)| \ge 1$. 
Note that $N(x,i)$ is a vertex cut of $G$ for each $1 \le i \le t-1$.
Then $|N(x,i)| \ge k$ for $1 \le i \le t-1$ since $G$ is $k$-connected. 
Therefore, \[n  \ge  1 + k(t-1) +1= k t - k+2.\] 
It follows that $ecc(x)=t \le \lfloor \frac{n+k-2}{k} \rfloor$.
Then \[\mbox{diam}(G)=\max\limits_{x \in V(G)}\{\mbox{ecc}(x)\} \le  \lfloor \frac{n+k-2}{k} \rfloor.\] 
This upper bound is sharp can be seen as follows. 
If $k=1$, then $\mbox{diam}(G) \le n-1$, and this upper bound can be obtained if $G$ is a path $P_{n}$ of order $n$.
If $k \ge 2$ is even, then by Lemma \ref{L:H(k,n)kEven},  this upper bound can be obtained
if $G$ is Harary graph $H_{k,n}$.  
In particular, if $k=2$, then $\mbox{diam} \le \lfloor \frac{n}{2} \rfloor$, 
and this upper bound can be obtained if $G$ is a cycle $C_{n}(=H_{2,n})$ of order $n$. \qed\\

\begin{theorem}\rm \label{T:W(x,G)}
Let $n-1>k \ge 1$ and $G$ be a $k$-connected graph of order $n$. 
Then for any vertex $x$ of $G$,
\[W(x,G) \le \frac{1}{2}\lfloor \frac{n+k-2}{k} \rfloor (2n+k-2-k\lfloor \frac{n+k-2}{k} \rfloor).\]
Moreover, this bound is sharp when $k=1$ or $k \ge 2$ is even. 
\end{theorem}

\proof Let $x$ be a vertex of $G$ with $\mbox{ecc}(x)=t$ where $t>1$.
Such a vertex exists since $G$ is not a complete graph.
Assume that $N(x,i)$ is the set of vertices with distance $i$ from $x$ in $G$ for $1 \le i \le t$. 
Then by Lemma \ref{L:W(x,G)N(x,i)},
\begin{eqnarray*}
W(x,G) &=& \sum\limits_{i=1}^{t} i \cdot |N(x,i)|= \sum\limits_{i=1}^{t} (t-(t-i)) \cdot |N(x,i)| \\
&=& t \sum\limits_{i=1}^{t} |N(x,i)| - \sum\limits_{i=1}^{t-1} (t-i) |N(x,i)|.
\end{eqnarray*}

Since $\sum\limits_{i=1}^{t} |N(x,i)|=n-1$ and $|N(x,i)| \ge k$ for $1 \le i \le t-1$,
\begin{eqnarray*}
W(x,G) &\le& (n-1) t - k \sum\limits_{i=1}^{t-1}(t-i)= (n-1) t - k \sum\limits_{i=1}^{t-1}i \\
&=& (n-1)t -k {t \choose 2} \\
&=& -\frac{k}{2} t^{2} +\frac{2n+k-2}{2} t.
\end{eqnarray*}

We distinguish three cases based on $k\ge 3$, $k=2$, and $k=1$ respectively.

Case 1. $k \ge 3$.   By Theorem \ref{T:D(G)}, 
$t \le \mbox{diam}(G)  \le  \lfloor \frac{n+k-2}{k} \rfloor$.

Note that $W(x,G) \le -\frac{k}{2} t^{2} +\frac{2n+k-2}{2} t$ is maximized at 
the integer which is closest to $\frac{2n+k-2}{2k}$ and at most $\lfloor \frac{n+k-2}{k} \rfloor$. 
Let $k\geq 3$. We want to show that an integer closest to 
\[\frac{2n+k-2}{2k}=\frac{n-2}{k}+\frac{1}{2}+\frac{1}{k}\] is 
\[t_{0}= \lfloor \frac{n+k-2}{k}\rfloor =\lfloor \frac{n-2}{k}+1 \rfloor \]

It is enough to show that
\[ \lfloor \frac{n-2}{k} + 1 \rfloor - \frac{1}{2} \leq 
\frac{n-2}{k}+\frac{1}{2}+\frac{1}{k} \leq \lfloor 
\frac{n-2}{k}+1\rfloor + \frac{1}{2}\]

The left side of the inequality follows because  $\lfloor \frac{n-2}{k} 
+1\rfloor \leq \frac{n-2}{k}+1+\frac{1}{k}$.  The right side of the 
inequality follows because $\frac{n-2}{k}<\lfloor \frac{n-2}{k} 
+1\rfloor$, hence $\frac{n-2}{k} +\frac{1}{k} \leq \lfloor \frac{n-2}{k} 
+1\rfloor$, since the smallest integer greater than $\frac{n-2}{k}$ can 
be reached by adding units of size $\frac{1}{k}$.

Bring the above integer $t_{0}= \lfloor \frac{n+k-2}{k}\rfloor $ into the upper bound on
$W(x,G)$, we have
\begin{eqnarray*}
W(x,G) & \le & -\frac{k}{2} t_{0}^{2} +\frac{2n+k-2}{2} t_{0}\\
&=&\frac{1}{2}\lfloor \frac{n+k-2}{k} \rfloor (2n+k-2-k\lfloor \frac{n+k-2}{k} \rfloor).
\end{eqnarray*}

By Lemma \ref{L:H(k,n)kEven},  this upper bound can be realized
by any vertex of Harary graph $G=H_{k,n}$ when $k$ is even.

Case 2. $k=2$. $W(x,G) \le -t^{2} + nt$.

As a function of $t$,  $-t^{2} + nt$ is quadratic and attains the maximum at $t_{m}=\frac{n}{2}$.
By Theorem \ref{T:D(G)}, $t \le \mbox{diam}(G) \le \lfloor\frac{n}{2}\rfloor \le t_{m}$.
It is easy to check that $W(x,G)$ attains its maximum $\lfloor \frac{n}{2}\rfloor \lceil \frac{n}{2}\rceil$ 
at $t_{0}=\lfloor\frac{n}{2}\rfloor=\lfloor \frac{n+k-2}{k} \rfloor$.
This upper bound can be realized by any vertex of $G=H_{2,n}$ (a cycle of order $n$).

Case 3. $k=1$.  $W(x,G) \le  -\frac{1}{2} t^{2} +\frac{2n-1}{2} t$.

As a function of $t$,  $-\frac{1}{2} t^{2} +\frac{2n-1}{2} t$ is quadratic and attains the maximum at $t_{m}=\frac{2n-1}{2}$. 
By Theorem \ref{T:D(G)}, $t \le \mbox{diam}(G) \le n-1 \le t_{m}$.
It is easy to check that $W(x,G)$ attains its maximum $\frac{1}{2} (n^{2}-n)$ 
at $t_{0}=n-1=\lfloor \frac{n+k-2}{k} \rfloor$.
This upper bound can be realized by an end vertex of $G=P_{n}$ (a path of order $n$).
\qed\\

\begin{theorem}\rm\label{T:W(G)}
Let $n-1>k \ge 1$ and $G$ be a $k$-connected graph of order $n$.
Then \[W(G) \le \frac{1}{4} n \lfloor \frac{n+k-2}{k} \rfloor (2n+k-2-k\lfloor \frac{n+k-2}{k} \rfloor).\]
Moreover, this bound is sharp when $k \ge 2$ is even. 
\end{theorem} 
\proof It follows immediately from Theorem \ref{T:W(x,G)}
and the fact that $W(G)=\frac{1}{2}\sum\limits_{x \in V(G)} W(x, G)$. 
Moreover, by Lemma \ref{L:H(k,n)kEven},  
this upper bound can be obtained when $G$ is Harary graph $H_{k,n}$, where $k \ge 2$ is even. 
\qed\\

\section{Final comments}

The authors made a final revision to this paper  in 2013.  
We did not publish our work as we realized that our results 
were stated equivalently under the terminology of mean distance by Favaron et al. 
as a brief Remark 2.6.1 \cite {FKM89} in 1989 without reference papers:
$\mu(G) \le \lfloor \frac{n+k-1}{k} \rfloor \frac{n-1-\frac{k}{2} \lfloor \frac{n-1}{k} \rfloor}{n-1}$ and this bound is attained by 
the  $\frac{k}{2}$th power of a cycle (that is, Harary graph $H_{k,n}$) when $k$ is even. 


\end{document}